\documentclass{amsart} 
\usepackage{url,graphicx}
\usepackage[dvips]{epsfig}
\usepackage{latexsym,color}
\usepackage{amscd,amssymb,verbatim}

\newcommand{\bd}{{\textrm{Bd}\,}}
\newcommand{\be}{\begin{enumerate}}

\newcommand{\bo}{\partial\,}
\newcommand{\bu}{\bullet}

\newcommand{\ca}{{\mathcal A}}

\newcommand{\cat}{\text{\bf Cat}}
\newcommand{\cb}{{\mathcal B}}

\newcommand{\cf}{{\mathcal F}}

\newcommand{\cm}{{\mathcal M}}

\newcommand{\co}{{\mathcal O}}

\newcommand{\cs}{{\mathcal S}}

\newcommand{\da}{\Delta}

\newcommand{\ee}{\end{enumerate}}

\newcommand{\grp}{\text{\bf Grp}}

\newcommand{\id}{\text{id}}
\newcommand{\il}{\text{colim}\,}

\newcommand{\lam}{\lambda}

\newcommand{\lr}{\leftrightarrow}

\newcommand{\mqed}{\square}

\newcommand{\pr}{\noindent{\bf Proof. }}

\newcommand{\ra}{\rightarrow}
\newcommand{\Ra}{\Rightarrow}

\newcommand{\sm}{\setminus}

\newcommand{\ssc}{\text{\bf SS}}

\newcommand{\tda}{\Delta}

\newcommand{\ti}{\tilde }
\newcommand{\tp}{\text{\bf Top}}

\newtheorem{thm}{Theorem}[section]
\newtheorem{df}  [thm]{Definition}

\newtheorem{prop}[thm]{Proposition}

\newtheorem{rem} [thm]{Remark}

\newtheorem{exam}[thm]{Example}
\numberwithin{equation}{section}

\title{Group Actions on Posets}

\author{Eric Babson and Dmitry N.\ Kozlov}

\date{\today \\[0.1cm]
  Research at MSRI was supported in part by NSF grant DMS-9022140.
  The first author was supported by an NSF Postdoctoral Fellowship.
  The second author was supported by the Swedish Science Council
  Postdoctoral Fellowship M-PD 11292-303}

\address{Department of Mathematics, University of Washington, Seattle,
  WA 98105, U.S.A.} \address{Department of Mathematics, Royal
  Institute of Technology, Stockholm, S-100 44, Sweden.}

\email{babson@@math.washington.edu, kozlov@@math.kth.se} 
 
\begin{document}

\begin{abstract}

In this paper we study quotients of posets by group actions. In
order to define the quotient correctly we enlarge the considered class
of categories from posets to loopfree categories: categories without
nontrivial automorphisms and in\-ver\-ses. We view group actions as
certain functors and define the quotients as co\-li\-mits of these
functors. The advantage of this definition over studying the
quo\-ti\-ent poset (which in our language is the colimit in the poset
category) is that the realization of the quotient loopfree category is
more often homeomorphic to the quotient of the realization of the
original poset. We give conditions under which the quotient commutes
with the nerve functor, as well as conditions which guarantee that the
quotient is again a~poset.

\end{abstract}

\maketitle

                 \section{Introduction}

Assume that we have a finite group $G$ acting on a poset $P$ in an
order-preserving way. The~purpose of this chapter is to compare the
various constructions of the quotient, associated with this action.
Our~basic suggestion is to view $P$ as a category and the group action
as a functor from $G$ to $\cat$. Then, it is natural to define $P/G$
to be the colimit of this functor.  As a~result $P/G$ is in general 
a~category, not a~poset.
                 
After getting a hand on the formal setting in Section~\ref{s5.2} we
proceed in Section~\ref{s5.3} with imposing different conditions on
the group action. We give conditions for each of the following
properties to be satisfied:
\begin{enumerate}
\item[(1)] the morphisms of $P/G$ are exactly the orbits of the
  morphisms of~$P$, we call it {\it regularity}; 
\item[(2)] the quotient construction commutes with Quillen's
  nerve functor;
\item[(3)]  $P/G$ is again a poset.
\end{enumerate}
  
Furthermore, we study the class of categories which can be seen as the
``quotient closure'' of the set of all finite posets: loopfree
categories.

In Section~4 we draw connections to determining the multiplicity of
the trivi\-al character in the induced representations of $G$ on the
homology groups of the nerve of the category, derive a formula for the
M\"obius function of $P/G$ and, based on formulae of Sundaram and
Welker,~\cite{SW}, give a quotient analog of Goresky-MacPherson
formulae.

As another example where these methods proved to be essential we would
like to mention the computation of the homology groups of the deleted
symmetric join of an~infinite simplex, see~\cite{BK}.

\section{Formalization of group actions and the main question}
\label{s5.2}
 
\subsection{Preliminaries} $\,$


\noindent
For a~small category $K$ denote the~set of its objects by $\co(K)$ and
the~set of its morphisms by $\cm(K)$. For every $a\in\co(K)$ there is
exactly one identity morphism which we denote $\id_a$, this allows us
to identify $\co(K)$ with a~subset of $\cm(K)$. If $m$ is a~morphism
of $K$ from $a$ to $b$, we write $m\in\cm_K(a,b)$, $\bo^\bu m=a$ and
$\bo_\bu m=b$. The~morphism $m$ has an~inverse $m^{-1}\in\cm_K(b,a)$,
if $m\circ m^{-1}=\id_a$ and $m^{-1}\circ m =\id_b$. If only
the~identity mor\-phisms have inverses in $K$ then $K$ is said to be
a~category without inverses.
  
We denote the~category of all small categories by $\cat$. If
$K_1,K_2\in\co(\cat)$ we denote by $\cf(K_1,K_2)$ the set of functors
from $K_1$ to $K_2$. We need three full subcategories of $\cat$: ${\bf
  P}$ the~category of posets, (which are categories with at most one
morphism, denoted $(x\ra y)$, between any two objects $x,y$), ${\bf
  L}$ the~category of loopfree categories (see Definition~\ref{5dc}),
and $\grp$ the~category of groups, (which are categories with a single
element, morphisms given by the group elements and the law of
composition given by group multiplication).  Finally, {\bf 1} is
the~terminal object of $\cat$, that is, the~category with one element,
and one (identity) morphism. The~other two categories we use are
$\tp$, the~category of topological spaces, and $\ssc$, the~category of
simplicial sets.

We are also interested in the~functors $\tda:\cat\ra\ssc$ and ${\mathcal
  R}:\ssc\ra\tp$.  The~composition is denoted $\ti\da:\cat\ra\tp$.
Here, $\da$ is the~nerve functor, see Appendix~B,
or~\cite{Q73,Q78,Se}. In particular, the~simplices of $\tda(K)$
are chains of morphisms in $K$, with dege\-ne\-rate simplices
corresponding to chains that include identity morphisms,
see~\cite{GM,W}. ${\mathcal R}$ is the~topological realization
functor, see~\cite{Mil}.

Note that both $\tda$ and $\ti\da$ have weak homotopy inverses,
i.e.,~functors $\xi:\ssc\ra\cat$ and $\ti\xi:\tp\ra\cat$ such that
$\ti\da\circ\ti\xi$ is homotopic to the~identity, and ${\mathcal
  R}\circ\tda\circ\xi$ is homotopic to ${\mathcal R}$, see~\cite{FL}.

We recall here the definition of a colimit (see~\cite{Ma,Mit}).  
\begin{df} \label{5col}
  Let $K_1$ and $K_2 $ be categories and
  $X\in\cf(K_1,K_2)$. A~{\bf sink} of $X$ is a~pair consisting of
  $L\in\co(K_2)$, and a~collection of morphisms $\{\lam_s\in\cm_{K_2}
  (X(s),L)\}_{s\in\co(K_1)}$, such that if
  $\alpha\in\cm_{K_1}(s_1,s_2)$ then $\lam_{s_2}\circ
  X(\alpha)=\lam_{s_1}$. (One way to think of this collection of
  morphisms is as a natural transformation between the functors $X$
  and $X'=X_1\circ X_2$, where $X_2$ is the terminal functor
  $X_2:K_1\ra\bf 1$ and $X_1:{\bf 1}\ra K_2$ takes the object of
  $\,\bf 1$ to $L$). When $(L,\{\lam_s\})$ is universal with respect
  to this property we call it the {\bf colimit} of $X$ and write
  $L=\il X$.
\end{df}

\subsection{Definition of the quotient and formulation of 
  the main problem} $\,$


\noindent
  Our main object of study is described in the~following definition.

\begin{df}\label{5adf}
  We say that a group $G$ {\bf acts on} a category $K$ if there is a
  functor $\ca_K:G\ra\cat$ which takes the~unique object of G to $K$.
  The~colimit of $\ca_K$ is called the {\bf quotient} of $K$ by the
  action of $G$ and is denoted by $K/G$.
\end{df}

To simplify notations, we identify $\ca_K g$ with $g$ itself.
Furthermore, in Definition~\ref{5adf} the category $\cat$ can be
replaced with any category $C$, then $K,K/G\in\co(C)$. Important
special case is $C=\ssc$. It arises when $K\in\co(\cat)$ and we
consider $\il\tda\circ\ca_K=\tda(K)/G$.

\vskip8pt
\noindent
{\bf Main Problem.} $\,$ {\it Understand the relation between the
  topological and the categorical quotients, that is, between
  $\tda(K/G)$ and $\tda(K)/G$.}
\vskip8pt

To start with, by the~universal property of colimits there exists
a~canonical sur\-jec\-tion $\lam:\tda(K)/G\ra\tda(K/G)$. In the~next
section we give com\-bi\-na\-torial con\-di\-tions under which this
map is an isomorphism.

The general theory tells us that if $G$ acts on the~category $K$, then
the~colimit $K/G$ exists, since $\cat$ is cocomplete. We shall now
give an explicit description.

\vskip4pt

\noindent
{\bf An explicit description of the category $K/G$.}

\vskip4pt

\noindent
When $x$ is a morphism of~$K$, denote by $Gx$ the~orbit of $x$ under
the~action of $G$. We have $\co(K/G)=\{Ga\,|\,a\in\co(K)\}$. The
situation with morphisms is more complicated.  Define a relation $\lr$
on the~set $\cm(K)$ by setting $x\lr y$, iff there are decompositions
$x=x_1\circ\dots\circ x_t$ and $y=y_1\circ\dots\circ y_t$ with
$Gy_i=Gx_i$ for all $i\in [t]$. The~relation $\lr$ is reflexive and
symmetric since $G$ has identity and inverses, however it is not in
general transitive.  Let $\sim$ be the~transitive closure of $\lr$, it
is clearly an equivalence relation. Denote the~$\sim$ equivalence
class of $x$ by $[x]$. Note that $\sim$ is the~minimal equivalence
relation on $\cm(K)$ closed under the $G$ action and under
composition; that is, with $a\sim ga$ for any $g\in G$, and if $x\sim
x'$ and $y\sim y'$ and $x\circ x'$ and $y\circ y'$ are defined then
$x\circ x'\sim y \circ y'$. It is not difficult to check that the~set
$\{[x]\,|\,x\in\cm(K)\}$ with the~relations $\bo_\bu[x]=[\bo_\bu x]$,
$\bo^\bu[x]=[\bo^\bu x]$ and $[x]\circ[y]=[x\circ y]$ (whenever
the~composition $x\circ y$ is defined), are the morphisms of
the~category $K/G$.

\vskip4pt

Note that if $P$ is a poset with a $G$ action, the~quotient taken in
$\cat$ need not be a poset, and hence may differ from the~poset
quotient.
\begin{exam} 
  Let $P$ be the center poset in the figure below. Let $\cs_2$ act on
  $P$ by simultaneously permuting $a$ with $b$ and $c$ with $d$. (I)
  shows $P/\cs_2$ in ${\bf P}$ and (II) shows $P/\cs_2$ in $\cat$.
  Note that in this case the quotient in $\cat$ commutes with the
  functor $\da$ (the canonical surjection $\lam$ is an isomorphism),
  whereas the quotient in ${\bf P}$ does not.
\end{exam}

$$\epsffile{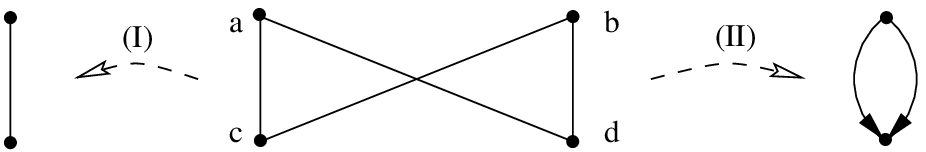}$$
$$\text{Figure 5.1}$$

  \section{Conditions on group actions} \label{s5.3}

\subsection{Outline of the results and surjectiveness of
the canonical map} $\,$
  

\noindent
In this section we consider combinatorial conditions for a~group $G$
acting on a~category $K$ which ensure that the quotient by the group
action commutes with the nerve functor. If $\ca_K:G\ra\cat$ is a~group
action on a~category $K$ then $\tda\circ\ca_K:G \ra \ssc$ is the
associated group action on the nerve of $K$.  It is clear that
$\tda(K/G)$ is a~sink for $\tda\circ\ca_K$, and hence, as previously
mentioned, the universal property of colimits gives a~canonical map 
$\lam:\tda(K)/G\ra\tda(K/G)$. We wish to find conditions under which 
$\lam$ is an isomorphism.
  
First we prove in Proposition~\ref{5sur} that $\lam$ is always
surjective. Furthermore, $Ga=[a]$ for $a\in\co(K)$, which means that,
restricted to $0$-skeleta, $\lam$ is an~iso\-morphism. If the two
simplicial spaces were simplicial complexes (only one face for any
fixed vertex set), this would suffice to show isomorphism. Neither one
is a~simplicial complex in general, but while the quotient of
a~complex $\tda(K)/G$ can have simplices with fairly arbitrary face
sets in common, $\tda(K/G)$ has only one face for any fixed edge set,
since it is a~nerve of a~category. Thus for $\lam$ to be
an~iso\-morphism it is necessary and sufficient to find conditions
under which


\be
\item[{1)}] $\lam$ is an isomorphism restricted to $1$-skeleta;
\item[{2)}] $\tda(K)/G$ has only one face with any given set of edges.
  \ee

We will give conditions equivalent to $\lambda$ being an~isomorphism,
and then give some stronger conditions that are often easier to check,
the strongest of which is also inherited by the action of any subgroup
$H$ of $G$ acting on $K$.

First note that a~simplex of $\tda(K/G)$ is a~sequence
$([m_1],\dots,[m_t])$, $m_i\in\cm(K)$, with $\bo_\bu [m_{i-1}]=\bo^\bu
[m_i]$, which we will call a~{\it chain}. On the other hand a simplex
of $\tda(K)/G$ is an~orbit of a~sequence $(n_1,\dots,n_t)$,
$n_i\in\cm(K)$, with $\bo_\bu n_{i-1}=\bo^\bu n_i$, which we denote
$G(n_1,\dots,n_t)$. The~canonical map $\lam$ is given by
$\lam(G(n_1,\dots,n_t))=([n_1],\dots,[n_t])$.

\begin{prop} \label{5sur}
  Let $K$ be a category and $G$ a group acting on $K$. The canonical map
$\lam:\tda(K)/G\ra\tda(K/G)$ is surjective.  
\end{prop}
\pr By the above description of $\lam$ it suffices to fix a chain
$([m_1],\dots,[m_t])$ and find a chain $(n_1,\dots, n_t)$ with
$[n_i]=[m_i]$. The proof is by induction on $t$. The case $t=1$
is obvious, just take $n_1=m_1$.

Assume now that we have found $n_1,\dots,n_{t-1}$, so that
$[n_i]=[m_i]$, for $i=1,\dots,t-1$, and $n_1,\dots,n_{t-1}$ compose,
i.e., $\bo^\bu n_i=\bo_\bu n_{i+1}$, for $i=1,\dots,t-2$. Since
$[\bo_\bu n_{t-1}]=[\bo_\bu m_{t-1}]=[\bo^\bu m_t]$, we can find $g\in
G$, such that $g\bo^\bu m_t=\bo_\bu n_{t-1}$. If we now take $n_t=g
m_t$, we see that $n_{t-1}$ and $n_t$ compose, and $[n_t]=[m_t]$,
which provides a proof for the induction step. 
\qed


\subsection{Conditions for injectiveness of the canonical projection} $\,$ 

\begin{df} \label{5df3.2}
  Let $K$ be a category and $G$ a group acting on $K$. We say that
  this action satisfies {\bf Condition (R)} if the following is true:
  If $x,y_a,y_b\in\cm(K)$, $\bo_\bu x=\bo^\bu y_a= \bo^\bu y_b$ and
  $Gy_a=Gy_b$, then $G(x\circ y_a)=G(x\circ y_b)$.  
\end{df}

$$\epsffile{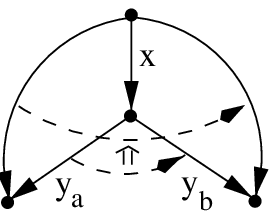}$$
$$\text{Figure 5.2}$$

We say in such case that $G$ acts {\it regularly} on~$K$.

\begin{prop} \label{5reg}
  Let $K$ be a category and $G$ a group acting on $K$. This action
  satisfies Condition (R) iff the~canonical surjection $\lam:\tda(K)/G
  \ra \tda(K/G)$ is injective on $1$-skeleta.
\end{prop}
\pr The injectiveness of $\lam$ on $1$-skeleta is equivalent to
requiring that $Gm=[m]$, for all $m\in\cm(K)$, while Condition~(R) is
equivalent to requiring that $G(m \circ Gn)=G(m \circ n)$, for all
$m,n\in\cm(K)$ with $\bo_\bu m=\bo^\bu n$; here $m\circ Gn$ means the
set of all $m\circ gn$ for which the composition is defined.

Assume that $\lam$ is injective on $1$-skeleta. The we have the
following computation:
$$G(m \circ Gn)=Gm \circ Gn=[m]\circ [n]=[m \circ n]=G(m \circ n),$$
hence the Condition (R) is satisfied.

Reversely, assume that the Condition~(R) is satisfied, that is $G(m
\circ Gn)=G(m \circ n)$.  Since the equivalence class $[m]$ is
generated by $G$ and composition, it suffices to show that orbits are
preserved by com\-po\-sition, which is precisely $G(m \circ Gn)=G(m
\circ n)$. \qed

\vskip4pt

The following theorem is the main result of this chapter. It provides
us with combinatorial conditions which are equivalent to $\lambda$
being an isomorphism.

\begin{thm} \label{5tc}
   Let $K$ be a category and $G$ a group acting on $K$. 
The~following two assertions are equivalent for any $t\geq 2$:
\begin{enumerate}
\item[(1$_t$)] {\bf Condition (C$_{\text{\bf t}}$).}  If
  $m_1,\dots,m_{t-1},m_{a},m_{b}\in\cm(K)$ with $\bo^\bu m_i=\bo_\bu
  m_{i-1}$ for all $2\leq i\leq t-1$, $\bo^\bu m_{a}=\bo^\bu
  m_{b}=\bo_\bu m_{t-1}$, and $Gm_{a}=Gm_{b}$, then there is some
  $g\in G$ such that $gm_{a}=m_{b}$ and $g m_i=m_i$ for $1 \leq i \leq
  t-1$.
\item[(2$_t$)] The~canonical surjection $\lam : \tda(K)/G \ra \tda(K/G)$ is 
injective on $t$-skeleta.  
\end{enumerate}
  In particular, $\lambda$ is an isomorphism iff {\bf (C$_{\text{\bf t}}$)} 
is satisfied for all $t\geq 2$. If this is the case, we say that Condition 
{\bf (C)} is satisfied
\end{thm}

$$\epsffile{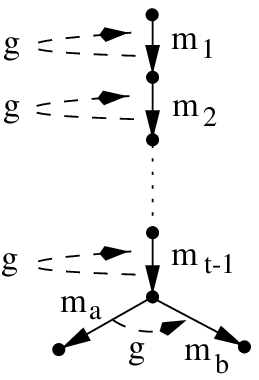}$$
$$\text{Figure 5.3}$$

\pr $(1_t)$ is equivalent to
$G(m_1,\dots,m_t)=G(m_1,\dots,m_{t-1},Gm_t)$; this notation is used,
as before, for all sequences $(m_1,\dots,m_{t-1},gm_t)$ which are
chains, that is for which $m_1\circ\dots\circ m_{t-1}\circ gm_t$ is
defined. $(2_t)$ implies Condition (R) above, and so can be restated
as $G(m_1,\dots,m_t)=(Gm_1,\dots,Gm_t)$.
\begin{multline*}
\underline{(2_t)\Ra (1_t):} \,\,G(m_1,\dots,m_t)=(Gm_1,\dots,Gm_t)=
G(Gm_1,\dots,Gm_t) \\
\supseteq G(m_1,\dots,m_{t-1},Gm_t)\supseteq G(m_1,\dots,m_t).  
\end{multline*}
\begin{multline*}
\underline{(1_2)\Ra (2_2):} \,\,G(m_1,m_2)=G(m_1,Gm_2) \\
=\{g_1(m_1,g_2 m_2)\,|\,\bo_\bu m_1=\bo^\bu g_2 m_2\}\, \\ 
=\{(g_1 m_1,g_2 m_2)\,|\,\bo_\bu g_1 m_1=\bo^\bu g_2 m_2\}=(Gm_1,Gm_2).
\end{multline*}

\hskip-2pt $\underline{(1_t)\Ra (2_t),\,t\geq 3:}$ We use induction on $t$.
\begin{equation*}
\begin{split}
 G(m_1, \dots, m_t)&=G(m_1, \dots, m_{t-1}, Gm_t) \\
&=\{(gm_1,\dots,gm_{t-1},\ti gm_t\,|\,\bo_\bu gm_{t-1}=\bo^\bu \ti gm_t)\} \\
&=\{(g_1 m_1,\dots,g_t m_t)\,|\,\bo_\bu g_i m_i=\bo^\bu g_{i+1}m_{i+1},
i\in[t-1]\} \\
&=(Gm_1,\dots,Gm_t).\,\,\,\,\mqed
\end{split}
\end{equation*}

\begin{exam}
  {\rm A~group action which satisfies Condition (C$_t$), but does not
    satisfy Condition~(C$_{t+1}$).}  Let $P_{t+1}$ be the order sum of
  $t+1$ copies of the 2-element antichain.  The~automorphism group of
  $P_{t+1}$ is the direct product of $t+1$ copies of ${\mathbb Z}_2$.
  Take $G$ to be the index $2$ subgroup consisting of elements with an
  even number of nonidentity terms in the product.
\end{exam}

\noindent  The~following condition implies Condition (C), and is often 
easier to check.

\vskip4pt

\noindent
{\bf Condition (S).} There exists a set $\{S_m\}_{m\in\cm(K)}$, 
$S_m\subseteq\,\,$Stab$\,(m)$, such that
\begin{enumerate}
\item[(1)] $S_m\subseteq S_{\bo^\bu m}\subseteq S_{m'}$, for any
  $m'\in\cm(K)$, such that $\bo_\bu m'=\bo^\bu m$;
\item[(2)] $S_{\bo^\bu m}$ acts transitively on 
$\{g m\,|\,g\in\,$Stab$\,(\bo^\bu m)\}$, for any $m\in\cm(K)$.
\end{enumerate}

\begin{prop} 
 Condition (S) implies Condition (C). 
\end{prop}
\pr Let $m_1,\dots,m_{t-1},m_a,m_b$ and $g$ be as in Condition (C),
then, since $g\in\,$Stab$\,(\bo^\bu m_a)$, there must exist $\ti g\in
S_{\bo^\bu m_a}$ such that $\ti g(m_a)=m_b$. From (1) above one can
conclude that $\ti g(m_i)=m_i$, for $i\in [t-1]$. \qed\vskip4pt

We say that the {\bf strong} Condition (S) is satisfied if Condition
(S) is satisfied with $S_a=\,$Stab$\,(a)$. Clearly, in such a~case
part (2) of the Condition~(S) is obsolete. 

\begin{exam}
  {\rm A~group action satisfying Condition~(S), but not the strong
    Condition~(S).}  Let $K=\cb_n$, lattice of all subsets of $[n]$
  ordered by inclusion, and let $G=\cs_n$ act on $\cb_n$ by permuting
  the ground set $[n]$. Clearly, for $A\subseteq[n]$, we have
  Stab$\,(A)=\cs_A\times\cs_{[n]\sm A}$, where, for $X\subseteq[n]$,
  $\cs_X$ denotes the subgroup of $\cs_n$ which fixes elements of
  $[n]\sm X$ and acts as a permutation group on the set $X$. Since
  $A>B$ means $A\supset B$, condition (1) of (S) is not satisfied for
  $S_A=\,$Stab$\,(A)$: $\cs_A\times\cs_{[n]\sm
    A}\not\supseteq\cs_B\times\cs_{[n]\sm B}$. However, we can set
  $S_A=\cs_A$. It is easy to check that for this choice of
  $\{S_A\}_{A\in\cb_n}$ Condition~(S) is satisfied.
\end{exam}

We close the discussion of the conditions stated above by the
following proposition.

\begin{prop}\label{5int} $\,$
  
\noindent  1) The~sets of group actions which satisfy Condition (C) or
  Condition (S) are closed under taking the restriction of the group
  action to a subcategory.

\noindent 2) Assume a finite group $G$ acts on a poset $P$, so that Condition
(S) is satisfied. Let $x\in P$ and $S_x\subseteq
H\subseteq\,$Stab$\,(x)$, then Condition (S) is satisfied for the
action of $H$ on $P_{\leq x}$.

\noindent 3) Assume a finite group $G$ acts on a category $K$, so that
Condition (S) is satisfied with $S_a=$Stab$(a)$ (strong version), and
$H$ is a subgroup of $G$.  Then the strong version of Condition $(S)$
is again satisfied for the action of $H$ on $K$.
\end{prop}
\pr 1) and 3) are obvious. To show 2) observe that for $a\leq x$ we
have $S_a\subseteq S_x\subseteq H$, hence $S_a\subseteq
H\cap\,$Stab$\,(a)$.  Thus condition (1) remains true. Condition~(2)
is true since
$\{g(b)\,|\,g\in\,$Stab$\,(a)\}\supseteq\{g(b)\,|\,g\in\,$Stab$\,(a)\cap
H\}$. \qed

\subsection{Conditions for the categories to be closed under taking
 quotients} $\,$


\noindent
Next, we are concerned with finding out what categories one may get as
a~quotient of a poset by a group action. In particular, we ask: {\it
  in which cases is the quotient again a poset?} To answer that
question, it is convenient to use the following class of
cate\-go\-ries.

\begin{df}\label{5dc}
  A~category is called {\bf loopfree} if it has no inverses and no
  nonidentity automorphisms.
\end{df}

Intuitively, one may think of loopfree categories as those which can
be drawn so that all nontrivial morphisms point down. To familiarize
us with the notion of a~loopfree category we make the following
observations:
\begin{itemize}
\item $K$ is loopfree iff for any $x,y\in\co(K)$, $x\neq y$, only one
  of the sets $\cm_K(x,y)$ and $\cm_k(y,x)$ is non-empty and
  $\cm_K(x,x)=\{\id_x\}$;
\item a poset is a loopfree category;
\item a barycentric subdivision of an arbitrary category is a loopfree
  category;
\item a barycentric subdivision of a loopfree category is a poset;
\item if $K$ is a loopfree category, then there exists a~partial order
  $\geq$ on the set $\co(K)$ such that $\cm_K(x,y)\neq\emptyset$
  implies $x\geq y$.
\end{itemize}

\begin{df}\label{5hor}
  Suppose $K$ is a~small category, and $T\in\cf(K,K)$. We say that $T$
  is {\bf hori\-zontal} if for any $x\in\co(K)$, if $T(x)\neq x$, then
  $\cm_K(x,T(x))=\cm_K(T(x),x)=\emptyset$. When a~group $G$ acts on
  $K$, we say that the action is hori\-zontal if each $g\in G$ is
  a~horizontal functor.
\end{df}

When $K$ is a~finite loopfree category, the action is always
horizontal. Another example of horizontal actions is given by rank
preserving action on a~(not necessarily finite) poset.  We have the
following useful property:

\begin{prop}\label{5star}
  Let $P$ be a finite loopfree category and $T\in\cf(P,P)$ be a
  hori\-zontal functor.  Let $\ti T\in\cf(\da(P),\da(P))$ be the
  induced functor, i.e., $\ti T=\da(T)$. Then $\da(P_T)=\da(P)_{\ti
    T}$, where $P_T$ denotes the subcategory of $P$ fixed by $T$ and
  $\da(P)_{\ti T}$ denotes the subcomplex of $\da(P)$ fixed by $\ti
  T$.
\end{prop}  

\pr Obviously, $\da(P_T)\subseteq\da(P)_{\ti T}$. On the other hand,
if for some $x\in\da(P)$ we have $\ti T(x)=x$, then the minimal
simplex $\sigma$, which contains $x$, is fixed as a set and, since the
order of simplices is preserved by $T$, $\sigma$ is fixed by $T$
pointwise, thus $x\in\da(P_T)$. \qed\vskip4pt

The~class of loopfree categories can be seen as the closure of the
class of posets under the operation of taking the quotient by a
horizontal group action.  More precisely, we have:
\begin{prop}
  The~quotient of a loopfree category by a~horizontal action is again
  a loopfree category. In particular, the quotient of a poset by
  a~horizontal action is a~loopfree category.
\end{prop}
\pr Let $K$ be a loopfree category and assume $G$ acts on $K$
horizontally. First observe that $\cm_{K/G}([x])=\{\id_{[x]}\}$.
Because if $m\in\cm_{K/G}([x])$, then there exist $x_1,x_2\in\co(K)$,
$\ti m\in\cm_K(x_1,x_2)$, such that $[x_1]=[x_2]$, $[\ti m]=m$. Then
$g x_1=x_2$ for some $g\in G$, hence, since $g$ is a horizontal
functor, $x_1=x_2$ and since $K$ is loopfree we get $\ti m=\id_{x_1}$.

Let us show that for $[x]\neq[y]$ at most one of the sets
$M_{K/G}([x],[y])$ and $M_{K/G}([y],[x])$ is nonempty. Assume the
contrary and pick $m_1\in M_{K/G}([x],[y])$, $m_2\in
M_{K/G}([y],[x])$. Then there exist $x_1,x_2,y_1,y_2\in\co(K)$, $\ti
m_1\in\cm_K(x_1,y_1)$, $\ti m_2\in\cm_K(y_2,x_2)$ such that
$[x_1]=[x_2]=[x]$, $[y_1]=[y_2]=[y]$, $[\ti m_1]=[m_1]$, $[\ti
m_2]=[m_2]$.  Choose $g\in G$ such that $g y_1=y_2$. Then $[g
x_1]=[x_2]=[x]$ and we have $g \ti m_1\in\cm_K(g x_1,y_2)$, so $\ti
m_2\circ g \ti m_1\in \cm_K(g x_1,x_2)$. Since $K$ is loopfree we
conclude that $g x_1=x_2$, but then both $\cm_K(x_2,y_2)$ and
$\cm_K(y_2,x_2)$ are nonempty, which contradicts to the fact that $K$
is loopfree. \qed\vskip4pt

Next, we shall state a~condition under which the quotient of
a~loopfree category is a~poset.

\begin{prop} \label{5sreg}
  Let $K$ be a loopfree category and let $G$ act on $K$. The~following
  two assertions are equivalent:
\begin{enumerate}
\item[(1)] {\bf Condition (SR).} If $x,y\in\cm(K)$, $\bo^\bu x=\bo^\bu
  y$ and $G\bo_\bu x=G\bo_\bu y$, then $Gx=Gy$.
\item[(2)] $G$ acts regularly on $K$ and $K/G$ is a poset.
\end{enumerate}
\end{prop}
\pr $(2)\Ra (1)$. Follows immediately from the regularity of the
action of $G$ and the fact that there must be only one morphism
between $[\bo^\bu x](=[\bo^\bu y])$ and $[\bo_\bu x](=[\bo_\bu y])$.

$$\epsffile{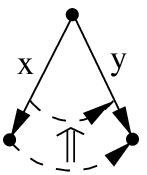}$$
$$\text{Figure 5.4}$$

$(1)\Ra (2)$. Obviously (SR) $\Ra$ (R), hence the action of $G$ is
regular. Furthermore, if $x,y\in\cm(K)$ and there exist $g_1,g_2\in G$
such that $g_1 \bo^\bu x=\bo^\bu y$ and $g_2 \bo_\bu x=\bo_\bu y$,
then we can replace $x$ by $g_1 x$ and reduce the situation to the one
described in Condition (SR), namely that $\bo^\bu x=\bo^\bu y$.
Applying Condition (SR) and acting with $g_1^{-1}$ yields the result. 
\qed

\vskip4pt

 When $K$ is a poset, Condition (SR) can be stated in simpler terms.

\vskip4pt

\noindent
{\bf Condition (SRP).} If $a,b,c\in K$, such that $a\geq b$, $a\geq c$
and there exists $g\in G$ such that $g(b)=c$, then there exists $\ti
g\in G$ such that $\ti g(a)=a$ and $\ti g(b)=c$. 

\[
\begin{array}{c}
\epsffile{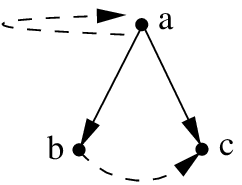}\\
\text{Figure 5.5}
\end{array}
\]


That is, for any $a,b\in P$, such that $a\geq b$, we require that
the stabilizor of $a$ acts transitively on $Gb$.

\begin{prop}
  Let $P$ be a poset and assume $G$ acts on $P$. The~action of $G$ on
  $P$ induces an action on the barycentric subdivision $\bd P$ (the
  poset of all chains of $P$ ordered by inclusion). This action
  satisfies Condition (S), hence it is regular and $\da(\bd
  P)/G\cong\da((\bd P)/G)$. Moreover, if the action of $G$ on $P$ is
  horizontal, then $(\bd P)/G$ is a poset.
\end{prop}

\pr Let us choose chains $b$, $c$ and $a=(a_1>\dots>a_t)$, such that
$a\geq b$ and $a\geq c$. Then $b=(a_{i_1}>\dots>a_{i_l})$,
$c=(a_{j_1}>\dots>a_{j_l})$. Assume also that there exists $g\in G$
such that $g(a_{i_s})=a_{j_s}$ for $s\in[l]$. If $g$ fixes $a$ then it
fixes every $a_i$, $i\in[t]$, hence $b=c$ and Condition (S) follows.

If, moreover, the action of $G$ is horizontal, then again
$a_{i_s}=a_{j_s}$, for $s\in[l]$, hence $b=c$ and Condition (SRP)
follows. \qed\vskip4pt

                 \section{Applications.}

  Let us first state two simple, but nevertheless fundamental, facts.

\begin{prop}\label{ap1}
  Assume that $G$ is a finite group which acts on $K$, a category
  without inverses, so that condition (C) is satisfied. Then
  $\beta_i(\da(K/G))=\langle\gamma_i,1\rangle$, where $\gamma_i$ is
  the induced representation of $G$ on $H_i(\da(K))$.
\end{prop}

\pr $\langle\gamma_i,1\rangle=\beta_i(\da(K)/G)=\beta_i(\da(K/G))$,
where the first equality follows from~\cite[Theorem 1]{Co} and the
second from Theorem~\ref{5tc}. \qed

\begin{prop}\label{**}
  Let $K$ be a finite loopfree category and $G$ a finite group which
  acts on~$K$, then
     $$\chi(\da(K)/G)=\dfrac{1}{|G|}\sum_{g\in G}\chi(\da(K_g)),$$
where $K_g$ denotes the subcategory of $K$ which is fixed by $g$.
\end{prop}

\pr $$\chi(\da(K)/G)=\dfrac{1}{|G|}\sum_{g\in G}\chi(\da(K)_g)=
\dfrac{1}{|G|}\sum_{g\in G}\chi(\da(K_g)),$$
where the first equality
follows from \cite[Theorem 2]{Co} and the second from
Proposition~\ref{5star}. \qed\vskip3pt

  Proposition~\ref{**} can be nicely restated in combinatorial language.
To do this we need the following definition.

\begin{df}\label{mudef}
  Let $K$ be a category, such that $\da(K)$ has finitely many
  simplices.  We define
  $\mu(K)\stackrel{\text{def}}{=}\ti\chi(\da(K))$. We call $\mu(K)$
  the {\bf M\"obius function} of $K$.
\end{df}

Clearly this definition generalizes the M\"obius function of a poset.
Similar definitions have been given: most notably (and apparently
independently) in~\cite{Hai} and~\cite{CLL}. We would like to mention
that if $K$ is a finite loopfree category then one has a
generalization of the recursive formula for the computation of the
M\"obius function (which is often taken as a definition of the
M\"obius function of a~poset):
$$\mu(\hat 0,x)=-\sum_{y\in\co(K),y<x}m_{x,y}\mu(\hat 0,y),$$
where
$m_{x,y}=|\cm_K(x,y)|$, $\mu(\hat 0,\hat 0)=1$ and $\mu(\hat 0,\hat
1)=\mu(K)$.  Here $\hat 0$ and $\hat 1$ are adjoint terminal and
initial objects.

\begin{prop}
  Let $K$ be a finite loopfree category, $G$ a finite group acting
  on~$K$, such that condition (C) is satisfied. Then
$$\mu(K/G)=\dfrac{1}{|G|}\sum_{g\in G}\mu(K_g).$$
\end{prop}
\pr Follows from Propositions~\ref{**} and~\ref{5tc} and the
definition of the M\"obius function. \qed\vskip3pt

As another application we obtain a quotient analog of the
Goresky-MacPherson formulae.

\begin{prop}\label{QGM}
  Let $\ca$ be a subspace arrangement in ${\mathbb C}^n$ and let $G$ be a
  finite group which acts on $\ca$. Assume that the induced action of
  $G$ on ${\mathcal L}_\ca$ (the intersection lattice of $\ca$) satisfies
  condition (C). Then
$$\beta^{n-1-i}(\cm_\ca/G)=\beta_i({\mathcal L}_\ca/G)=\sum_{x\in{\mathcal L}_\ca^{>\hat 0}/G}
\beta_{i-\dim x-1}(\da((\hat 0,x)/\,\text{Stab}\,(x))).$$
\end{prop}
\pr Follows from~\cite[Corollaries 2.8 and 2.10]{SW},
Theorem~\ref{5tc} and Proposition~\ref{5int}. \qed

\begin{rem}
  A~similar (though less pretty) formula can be derived for real
  subspace arrangements, cf.\ ~\cite[Theorems 2.4 and 2.5]{SW}.
\end{rem}

 \vskip3pt
\noindent
{\bf Acknowledgments.} We would like to thank Eva-Maria Feichtner for the 
careful reading of this paper.

\end{document}